\documentclass[a4paper,9pt]{article}
\usepackage{amssymb,amsmath,amsthm,graphicx,fancyhdr}

\newtheorem{Def}{Definition}
\newtheorem{Thm}{Theorem}

\newtheorem{Lem}{Lemma}[section]

\theoremstyle{remark}
\newtheorem{Rem}{Remark}[section]

\newcommand{\C}{\mathbb{C}}
\newcommand{\R}{\mathbb{R}}

\newcommand{\N}{\mathbb{N}}

\newcommand{\F}{\mathcal{F}}

\renewcommand{\Re}{\operatorname{Re}}
\newcommand{\I}{\infty}

\newcommand{\norm}[1]{\left\lVert #1\right\rVert}
\newcommand{\Lebn}[2]{\left\lVert #1 \right\rVert_{L^{#2}}}

\newcommand{\Jbr}[1]{\left\langle #1 \right\rangle}

\newcommand{\IN}{\quad\text{in }}

\def\({\left(}
\def\){\right)}
\def\<{\left\langle}
\def\>{\right\rangle}
\def\le{\leqslant}
\def\ge{\geqslant}

\def\d{{\partial}}
\def\l{\lambda}
\newcommand{\al}{\alpha}

\newcommand{\g}{\gamma}

\newcommand{\ka}{\kappa}

\newcommand{\eps}{\varepsilon}
\newcommand{\om}{\omega}

\newcommand{\la}{\lambda}

\newcommand{\de}{\delta}
\newcommand{\De}{\Delta}

\newcommand{\Om}{\Omega}

\title
{An example of stable excited state on nonlinear Schr\"odinger equation
with nonlocal nonlinearity}
\author{Masaya Maeda${}^{\dagger}$ and Satoshi Masaki${}^{\ddagger}$}

\date{}

\begin{document}

\maketitle

\vskip-5mm
\centerline{${}^{\dagger}$Mathematical institute, Tohoku University, 
}
\centerline{Sendai, 980-8578, Japan}

\vskip3mm

\centerline{${}^{\ddagger}$Department of Mathematics, 
Gakushuin University,}
\centerline{Toshima-ku Tokyo, 171-8588, Japan}

\begin{abstract}
In this article, we consider nonlinear Schr\"odinger equation with nonlocal nonlinearity 
which is a generalized model of the Schr\"odinger-Poisson system
(Schr\"odinger-Newton equations) in low dimensions.
We first prove the global well-posedness in wider space than in previous result and
show the stability of standing waves including excites states.
It turns out that an example of stable excite states with high Morse index is contained.
Several examples of traveling-wave-type solutions are also given.
\end{abstract}

 \section{Introduction}
In this article,
we consider the following nonlinear Schr\"odinger equation with nonlocal nonlinearity:
\begin{equation}\label{eq:nh}\tag{H}
	\left\{
	\begin{aligned}
	&2iu_t=-\De u +\la |x|^2u+\eta \(|x|^2*|u|^2\)u,\ (t,x)\in\R^{1+d}, \\
	& u(0,x)= u_0(x)\in \Sigma,
	\end{aligned}
	\right.
\end{equation}
where $\la,\eta\in\R$ and 
$\Sigma:=\{\phi\in H^1\ |\ \Jbr{x}\phi\in L^2\}$.
The operator $*$ stands for the usual convolution: $f*g(x)=\int_{\R^d} f(x-y)g(y)dy$.
It is shown in \cite{MHpre} that \eqref{eq:nh} is globally well-posed in $\Sigma$ and
further all the solution is written explicitly.
In this article, we give several examples of standing waves
and investigate stability of them.
This contains a new type example of a stable excited state.
\smallbreak
The equation \eqref{eq:nh} is a generalized model of 
the Schr\"odinger-Poisson system (or the Schr\"odinger-Newton equations) 
\begin{equation}\label{eq:sp}\tag{SP}
	\left\{
	\begin{aligned}
	&2iu_t=-\De u + Vu,\ (t,x)\in\R^{1+d}, \\
	&-\De V = |u|^2,\ x\in \R^d,\\
	& u(0,x)= u_0(x)
	\end{aligned}
	\right.
\end{equation}
in low dimensions.
Indeed, one sees that $V=-(1/2)(|x|*|u|^2)$ when $d=1$.
Generalizing with respect to growth order of interaction potential
and dimension, we reach to
\begin{equation}\label{eq:nhgen}
	\left\{
	\begin{aligned}
	&2iu_t=-\De u + \eta (|x|^\g*|u|^2) u,\ (t,x)\in\R^{1+d}, \\
	& u(0,x)= u_0(x).
	\end{aligned}
	\right.
\end{equation}
The equation \eqref{eq:nh} (with $\la=0$) corresponds to the case $\g=2$.
For the well-posedness result on \eqref{eq:sp},
see \cite{DLR07} for $d=1$ and see \cite{M2DSPpre} for $d=2$.
Well-posedness result of \eqref{eq:nhgen} for $0<\g\le 2$
is shown in \cite{MHpre}.
The feature of these equations is that the nonlinearity grows at the spatial infinity.
Further, with these results, it has turned out that
this kind of growing nonlinearity has an effect like a linear potential.
More precisely, the nonlinear potential $(W*|u|^2)(x)$ behaves like a linear potential
$W(x)\Lebn{u_0}2^2$ near the spatial infinity.
Hence, we put a linear potential $\l |x|^2$ in \eqref{eq:nh} in order to investigate
the competition between this potential and the one created
by the nonlinearity.

An idea to treat this kind of nonlinearity is the following:
We split the nonlinear potential of \eqref{eq:nh} as
\[
	(|x|^2*|u|^2)=|x|^2 \int |u(y)|^2dy  - 2x\cdot \int y|u(y)|^2 dy + \int|y|^2|u(y)|^2 dy.
\]
The point is that the first two terms are
unbounded in $x$ but time behavior of them can be specified thanks to 
the conservation laws. The last term can be absorbed by gauge transformation.
What requires the full
regularity assumption $u\in \Sigma$ is just the third term.
It therefore can be expected that if this part is removed then
the regularity assumption can be loosen.
From this respect, let us introduce
\begin{equation}\label{eq:nh2}\tag{H${}^\prime$}
	\left\{
	\begin{aligned}
	&2iv_t=-\De v +\la |x|^2v+\eta v\int_{\R^d} (|x-y|^2-|y|^2)|v(y)|^2dy
	,\ (t,x)\in\R^{1+d}, \\
	& v(0,x)= v_0(x)\in \Sigma^{1/2},
	\end{aligned}
	\right.
\end{equation}
where $\Sigma^{1/2}:=\{ \phi\in H^{1/2}\ |\ \Jbr{x}^{1/2}\phi \in L^2 \}$.
Notice that the modified nonlinearity makes sense
for $v\in \Sigma^{1/2}$.
When $v_0\in \Sigma^{1/2} \setminus \Sigma$ it turns out that the corresponding solution
has infinite energy.
Existence of infinite energy solution for 2D Schr\"odinger-Poisson system
is shown in \cite{MCPDE} with the same modification of nonlinear potential.

In this article, we turn to the study of standing waves of to \eqref{eq:nh} or \eqref{eq:nh2}.
Here, standing waves are solutions of the form $e^{-i\om t/2}\phi(x)$.
The number $\om\in\R$ is referred to as the frequency.
Standing waves of 1D Schr\"odinger-Poisson system are studied in \cite{CS07,DL10}.
The meaning ``a standing wave 
 is stable'', will be clear later, after we introduce the well-posedness results of \eqref{eq:nh} and \eqref{eq:nh2}.
For the stability of standing waves of general Hamiltonian PDEs, there is an elegant theory by Grillakis, Shatah, and Strauss \cite{GSS, GSS2}.
Let $E$ be a $C^2$ functional which satisfies $E(e^{is}u)=E(u)$ for all $s\in\R$,
and consider the equation
\begin{eqnarray}\label{eq:abs}
2iu_t=E'(u).
\end{eqnarray}
Set $Q(u)=\frac{1}{2}||u||_2^2$ and $S_{\om}=E-\om Q$, where $S_{\om}$ is called the ``action''.
Then, $S_{\om}'(\phi_{\om})=0$ holds if and only if $e^{-i\om t/2}\phi_{\om}$ is a standing wave solution.
We assume that there exists a frequency-to-standing-wave map $\om\mapsto \phi_{\om}$ such that $S_{\om}'(\phi_{\om})=0$ and $\mathrm{Ker}S_{\om}''(\phi_{\om})$ is spanned by $i\phi_{\om}$.
Let $n(S_{\om}''(\phi_{\om}))$ be the number of negative eigenvalues of $S_{\om}''(\phi_{\om})$.
What
Grillakis, Shatah, and Strauss \cite{GSS, GSS2} have shown are summarized as follows
\begin{enumerate}
\item[$(\mathrm{i})$]
If $n(S_{\om}''(\phi_{\om}))=0$, then $e^{-i\om t/2}\phi_{\om}$ is stable.
\item[$(\mathrm{ii})$]
If $n(S_{\om}''(\phi_{\om}))=1$ and $d''(\om)>0$, then $e^{-i\om t/2}\phi_{\om}$ is stable.
\item[$(\mathrm{iii})$]
If $d''(\om)\neq 0$ and $n(S_{\om}''(\phi_{\om}))-\max(d''(\om)/|d''(\om)|, 0)$ is odd, then $e^{-i\om t/2}\phi_{\om}$ is unstable.
\end{enumerate}
Obviously, there exist two more cases
\begin{enumerate}
\item[$(\mathrm{iv})$]
$n(S_{\om}''(\phi_{\om}))\geq 2$, $d''(\om)\neq 0$ and $n(S_{\om}''(\phi_{\om}))-\max(d''(\om)/|d''(\om)|, 0)$ is even.
\item[$(\mathrm{v})$]
The case $d''(\om)=0$.
\end{enumerate}
These cases are not discussed in \cite{GSS, GSS2}.
For the case (v), there are several results (\cite{CP03}, \cite{MDpre}, \cite{Ohta11}).
In the case (iv), the standing wave has high Morse indices.
So, $\phi_{\om}$ is not a minimizer of $E$ under the constraint of $Q$.
In this case, we say $\phi_{\om}$ is unstable in energy.
Standing waves which are unstable in energy seem to be unstable,
and in many cases they are actually unstable.
We shall show that, however, 
there exist stable standing waves belonging to the case (iv)
with arbitrary $n(S_{\om}''(\phi_{\om}))\geq 2$
(Theorem \ref{thm:main}). 

Our stability results are not based on the general theory.
We fully use a representation of the solution.
This is why we are able to discuss stability of excited states, which is hard to study in general.
Further, we can discuss the stability in a class in which
the energy functional does not necessarily take a finite value.

The results will be stated in the following order.
Global well-posedness results of \eqref{eq:nh} and \eqref{eq:nh2} are
given in Theorems \ref{thm:existence} and \ref{thm:existence2}, respectively.
We then discuss stability of standing waves of \eqref{eq:nh} and \eqref{eq:nh2}
in Theorem \ref{stability}.
In Theorem \ref{thm:main}, an example of excited state with high Morse index is introduced.
At last, in Theorems \ref{thm:single} and \ref{thm:multi},
we consider some more example of standing waves including
traveling-wave type solutions
$$e^{-i\om t/2}e^{ig_1(t)}e^{ix\cdot g_2(t)}\phi(x-g_3(t)),$$
and superpositions of finite number of such traveling waves.

\subsection{Existence and representation of a solution}
Before stability results,
we summarize well-posedness results for \eqref{eq:nh} and \eqref{eq:nh2}.
Set $M[u]=\Lebn{u}2^2$,
\begin{align*}
	X[u]&{}:=\int_{\R^d}x|u|^2\,dx, &
	P[u]&{}:=\mathrm{Im}\int_{\R^d}\bar{u}\nabla u\,dx=\int_{\R^d}\xi|\F{u}|^2\,d\xi.
\end{align*}
We also introduce energy
\[
	E[u] = \frac12 \Lebn{\nabla u}2^2 + \frac{\l}2 \Lebn{xu}2^2
	+\frac{\eta}4\iint_{\R^{d+d}} |x-y|^2 |u(x)|^2 |u(y)|^2 dx dy.
\]
Set ${\cal U}_{\ka}(t)=e^{i\frac{t}2(\Delta-\ka|x|^2)}$
Define an $\R^d$-valued function $g_\kappa(t)=g_\kappa(t,{\bf a},{\bf b})$
as a solution of an ODE $g_\kappa^{\prime\prime}=-\kappa g_\kappa$,
$g_\kappa^\prime(0)={\bf a}$, and $g_\kappa(0)={\bf b}$,
where $\kappa\in\R$ is a parameter.
Introduce the Galilean transform ${\cal G}_{\ka}(t)={\cal G}_{\ka}(t,{\bf a},{\bf b})$ by
\begin{equation}\label{eq:GT}
	({\cal G}_{\ka}(t,{\bf a},{\bf b})\phi)(x)
	= e^{-\frac{i}2 g_\ka(t,{\bf a},{\bf b})\cdot g_\ka^\prime(t,{\bf a},{\bf b})}
	e^{ix\cdot g_\ka^\prime(t,{\bf a},{\bf b})}
	\phi(x-g_\ka(t,{\bf a},{\bf b})).
\end{equation}
Basic properties of Galilean transform are summarized in Section \ref{sec:galilean}.

\begin{Thm}[\cite{MHpre}, Theorem 2.1]\label{thm:existence}
\eqref{eq:nh} is globally well-posed in $\Sigma$.
Uniqueness holds unconditionally.
Further, let $M=||u_0||_{2}^2$, ${\bf a}=P[u_0]/M$, ${\bf b}=X[u_0]/M$, and $\ka=\la+\eta M$.
Mass $M[u(t)]$ and energy $E[u(t)]$ are conserved.
Then the solution is written as
\begin{equation}\label{eq:mf}
\begin{aligned}
	u(t)={}&e^{-i\Psi(t)}{\cal G}_{\la}(t,{\bf a},{\bf b}){\cal U}_{\ka}(t){\cal G}_{\ka}(0,{\bf a},{\bf b})^{-1}u_0, \\
	={}& e^{-i\Psi(t)}{\cal G}_{\la}(t,{\bf a},{\bf b}){\cal G}_{\ka}(t,{\bf a},{\bf b})^{-1}{\cal U}_{\ka}(t) u_0,
\end{aligned}
\end{equation}
where $\Psi(t)$ is given by
$
	\Psi(t)=\frac{\eta}{2}\int_0^t||x{\cal U}_{\ka}(s){\cal G}_{\ka}(0,{\bf a},{\bf b})^{-1}u_0||_2^2\,ds.
$
\end{Thm}
\begin{Rem}
One easily verifies that $\Lebn{u(t)}p=\Lebn{\mathcal{U}_\kappa(t)u_0}p$
for all $p$
and that $X[u(t)]=X[\mathcal{U}_\la(t)u_0]$.
This means that dispersive properties of the solution is indicated by $\kappa$,
while the motion of the whole system is by $\la$.
\end{Rem}
\begin{Rem}
Unconditional uniqueness means that
the solution is unique in the class $C(\R:\Sigma)$ (without any additional assumption).
\end{Rem}
\begin{Thm}\label{thm:existence2}
\eqref{eq:nh2} is globally well-posed in $\Sigma^{1/2}$.
Uniqueness holds unconditionally.
Further, let $M=||v_0||_{2}^2$, ${\bf a}=P[v_0]/M$, ${\bf b}=X[v_0]/M$, and $\ka=\la+\eta M$.
Mass $M[u(t)]$ is conserved. Energy $E[u(t)]$ is finite and conserved if $v_0\in\Sigma$.
Then the solution is written as
\begin{equation}\label{eq:mf2}
\begin{aligned}
	v(t)={}& e^{-i\Phi(t)} {\cal G}_{\la}(t,{\bf a},{\bf b}){\cal U}_{\ka}(t){\cal G}_{\ka}(0,{\bf a},{\bf b})^{-1}v_0, \\
	={}& e^{-i\Phi(t)}{\cal G}_{\la}(t,{\bf a},{\bf b}){\cal G}_{\ka}(t,{\bf a},{\bf b})^{-1}{\cal U}_{\ka}(t) v_0,
\end{aligned}
\end{equation}
where $\Phi(t)=-\frac{\eta M}2 \int_0^t |g_\la(s,{\bf a},{\bf b})|^2 ds$.
\end{Thm}
\begin{Rem}
Theorem \ref{thm:existence2} is a generalization of
Theorem \ref{thm:existence} in such a sense that
$u_0=v_0\in \Sigma$ implies $u(t)=e^{i(\Psi(t)-\Phi(t))}v(t)$.
This can be seen from $\Lebn{xu}2^2=(2/\eta)(\Psi^\prime-\Phi^\prime)$.
Remark that if $v_0\in \Sigma^{1/2}\setminus \Sigma$ then
the energy of the corresponding solution is infinite.
\end{Rem}

\subsection{Stable excited state}
Here we state the results for stability of standing waves including excited states.
Let us first give a definition of the stability of standing waves.
\begin{Def}
Let $e^{-i\om t/2}\phi(x)$ be a solution of $(\ref{eq:nh})$.
We say $e^{-i\om t/2}\phi(x)$ is stable in a Banach space $X$
 if for all $\eps>0$, there exists $\de>0$ which satisfies the following.
If $||u_0-\phi||_{X}<\de$, then
\begin{eqnarray*}
\sup_{t>0}\inf_{s\in\R,y\in\R^d}||\phi-e^{is}u(t,\cdot -y)||_{X}<\eps,
\end{eqnarray*}
where $u(t)\in X$ is the solution of $\eqref{eq:nh}$
 with $u(0)=u_0\in X$.
\end{Def}
Stability of standing waves of \eqref{eq:nh2} is also defined by replacing
\eqref{eq:nh} and a solution $u$ in the above definition with 
\eqref{eq:nh2} and a solution $v$, respectively.
Notice that the definition of stability implicitly requires
well-posedness of the Cauchy problem.
In our case, they are already established in Theorems \ref{thm:existence}
and \ref{thm:existence2}. 

For a multi-index $\underline{n}=(n_1,\cdots,n_d)$, we introduce 
\begin{equation}
	\Omega_{\underline{n}} (x) := 
	\frac{1}{\sqrt{\pi^{{d}/{2}}
	2^{|\underline{n}|} \underline{n}!}}e^{\frac{|x|^2}2} 
	 \(-\frac{\d}{\d{x}}\)^{\underline{n}} e^{-|x|^2},
\end{equation}
where $\underline{n}!=\prod_{j=1}^d n_j!$.
Also set $\Omega_{\underline{n},\ka}(x)=\ka^{d/8} \Omega_{\underline{n}}(\ka^{1/4}x)$ for $\kappa>0$.
It is well known that, for any fixed $\ka>0$,
$\{\Omega_{\underline{n},\ka} \}_{\underline{n}}$ is a CONS of $L^2(\R^d)$,
and each $\Omega_{\underline{n},\ka}$ is an eigenfunction of $-\frac12\Delta
+\frac{\kappa}2|x|^2$ associated with an eigenvalue $\ka^\frac12 (|\underline{n}|+\frac{d}2)$.
For $s>0$, set $\Sigma^{s}:=D((-\frac12\Delta+\frac12|x|^2)^{s/2})$
and define $\Sigma^{-s}$ as a dual of $\Sigma^{s}$.
Set $\Sigma^0=L^2$.
Norm of $\Sigma^s$ ($s\in \R$) is given by
\[
	\norm{f}_{\Sigma^s} :=
	\(\sum_{\underline{n}} \(|\underline{n}|+ \frac{d}2 \)^{\frac{s}2}
	|(f,\Omega_{\underline{n}})|^2\)^{\frac12}.
\]
The space $\Sigma$ is identical to $\Sigma^1$.
Our stability result is the following.

\begin{Thm}\label{stability}
Let $d\geq 1$. Suppose $\la\geq0$, $\eta\in\R$ and $M>0$ satisfy $\kappa=\la +\eta M>0$.
For any multi-index $\underline{n}$, the following hold.
\begin{enumerate}
\item Set $\om_1=(1 + \eta M\ka^{-1})\ka^{\frac12}(|\underline{n}|+\frac{d}2)$.
Then, $e^{-i\om_1 t/2} M^{1/2}\Om_{\underline{n},\ka}(x)$ is a standing wave solution 
of \eqref{eq:nh} and is stable in $\Sigma^s$ for $s\geq 1$.
\item Set $\om_2=\ka^{\frac12}(|\underline{n}|+\frac{d}2)$.
Then, $e^{-i\om_2 t/2} M^{1/2}\Om_{\underline{n},\ka}(x)$ is a standing wave solution 
of \eqref{eq:nh2} and is stable in $\Sigma^{s}$ for $s\geq \frac{1}{2}$.
\end{enumerate}
\end{Thm}
One sees that the above standing wave solution 
is Ground state when $\underline{n}=0$, and is an excited state otherwise
(see Lemma \ref{lem:energy}).

Next theorem is concerned with stable excited states with high Morse indices.
Since our purpose is to show that such a state does exist,
let us restrict ourselves to the simplest situation.
So, we only treat the one dimensional case of \eqref{eq:nh},
and the case where a perturbation is supposed to be an even function.

\begin{Thm}\label{thm:main}
Let $d=1$ and set $L^2_r=\{ u\in L^2\ |\ u(x)=u(-x), x\in\R\}$.
\begin{enumerate}
\item[$(\mathrm{i})$]
Let $\la=0$ and $\eta=1$.
Then for any $m\in\N$, there exists a stable standing wave $e^{-i\om t/2}\phi_{\om}$ with $\phi_{\om}\in\Sigma_r$, $d''(\om)<0$ and $n(\left.S_{\om}''(\phi_{\om})\right|_{L^2_r})=2m$.
\item[$(\mathrm{ii})$]
Let $\la=2$ and $\eta=-1$.
Then for any $m\in\N$, there exists a stable standing wave $e^{-i\om t/2}\phi_{\om}$ with $d''(\om)>0$ and $n(\left.S_{\om}''(\phi_{\om})\right|_{L^2_r})=2m+1$.
\end{enumerate}
\end{Thm}

\subsection{Examples of standing waves}
Let us show some more example of standing wave solutions.
The first one is a traveling-wave-type solution with one peak.
\begin{Thm}[Single-peak standing wave]\label{thm:single}
Let $d\ge1$.
Suppose $\la,\eta \in \R$ and $M>0$ satisfy $\kappa=\la +\eta M>0$.
Suppose
$
	u_0 (x) =v_0(x)= M^\frac12 \mathcal{G}_\ka(0,{\bf a}_1, {\bf b}_1)
	\Omega_{\underline{n},\ka}(x),
$
where $\underline{n}$ is a multi-index and
${\bf a}_1, {\bf b}_1 \in \R^d$
\begin{enumerate}
\item The corresponding
solutions of \eqref{eq:nh} becomes
\begin{equation}\label{eq:single}
	u(t,x) = e^{-it\omega_1/2 }
	M^\frac12 \mathcal{G}_\la(t,{\bf a}_1, {\bf b}_1)
	\Omega_{\underline{n},\ka}(x)
\end{equation}
with $\omega_1=(\ka^\frac12 
+\eta M \ka^{-\frac12})(|\underline{n}|+\frac{d}2)$.
\item If $|{\bf a}_1|^2=\la |{\bf b}_1|^2$ then the corresponding
solutions of \eqref{eq:nh2} becomes
\begin{equation}\label{eq:single2}
	v(t,x) = e^{-it\omega_2/2 }
	M^\frac12 \mathcal{G}_\la(t,{\bf a}_1, {\bf b}_1)
	\Omega_{\underline{n},\ka}(x),
\end{equation}
with $\omega_2=(\ka^\frac12 
-\eta M |{\bf b}_1|^2)(|\underline{n}|+\frac{d}2)$.
\end{enumerate}
\end{Thm}

The next example is a solution which is a 
superposition of several traveling waves.
\begin{Thm}[Multi-peak standing wave]\label{thm:multi}
Let $d\ge1$.
Suppose $\la,\eta \in \R$ and $M>0$ satisfy $\kappa=\la +\eta M>0$.
Let $L$ be a positive integer.
For $\alpha_j \in \C$, ${\bf a}_j, {\bf b}_j \in \R^d$, and
$\underline{n}_j \in (\N\cup\{0\})^d$ ($j=1,2,\dots,L$), set
\[
	u_0 (x)=v_0(x) = \mu \sum_{j=1}^L \alpha_j \mathcal{G}_\ka(0,{\bf a}_j, {\bf b}_j)
	\Omega_{\underline{n}_j,\ka}(x),
\]
where $\mu>0$ is chosen so that $M[u_0]=M$.
Then, the solutions of \eqref{eq:nh} and \eqref{eq:nh2} are
\begin{align}\label{eq:multi}
	u(t,x) ={}& \mu \sum_{j=1}^L 
	\widetilde{\alpha}_j e^{-i \Psi_j(t)}
	\mathcal{G}_\la(t,{\bf a}, {\bf b})
	\mathcal{G}_\ka(t,{\bf a}_j-{\bf a}, {\bf b}_j-{\bf b})
	\Omega_{\underline{n}_j,\ka}(x) \\ \label{eq:multi2}
	v(t,x) ={}& \mu \sum_{j=1}^L 
	\widetilde{\alpha}_j e^{-i \Phi_j(t)}
	\mathcal{G}_\la(t,{\bf a}, {\bf b})
	\mathcal{G}_\ka(t,{\bf a}_j-{\bf a}, {\bf b}_j-{\bf b})
	\Omega_{\underline{n}_j,\ka}(x),
\end{align}
respectively, where ${\bf a}=P[u_0]/M$, ${\bf b}=X[u_0]/M$, 
$\widetilde{\alpha}_j=\alpha_j e^{\frac{i}2({\bf a}_j\cdot {\bf b}-{\bf b}_j\cdot {\bf a})}$,
$\Psi_j(t) = \Psi(t) + \frac{\ka^\frac12}{2}(|\underline{n}_j|+\frac{d}2)t$, and 
$\Phi_j(t) = \Phi(t) + \frac{\ka^\frac12}{2}(|\underline{n}_j|+\frac{d}2)t$.
\end{Thm}

\begin{Rem}
One sees from \eqref{eq:multi} that
each peak keeps its shape for all time.
The motion of a peak in physical and Fourier spaces is
described by a sum of two Galilean transform.
The first one is $\mathcal{G}_\ka(t,{\bf a}_j-{\bf a}, {\bf b}_j-{\bf b})$.
Recall that ${\bf a}$ and ${\bf b}$ denote
mean momentum and mean position of the whole system,
and so ${\bf a}_j-{\bf a}$ and ${\bf b}_j-{\bf b}$ are
relative momentum and relative position of each peak, respectively.
This transform therefore suggests that
all peaks are rotating around the center of mass
with the same period $2\pi \ka^{-1/2}$.
The second transform $\mathcal{G}_\la(t,{\bf a}, {\bf b})$
corresponds to the motion of the center of mass.
More precisely, $X[u(t)]=Mg_\la(t,{\bf a},{\bf b})$.
If $\la=0$ (no external potential) then the motion is a straight line, and
if $\la>0$ then the center of mass is rotation around the origin 
with period $2\pi \la^{-1/2}$.
Notice that, as long as $\eta\neq0$ (in presence of nonlinearity),
the period in which each peak rotates around the center of mass
and the period in which the center of mass does around the origin
do not coincide.
\end{Rem}

The rest of paper is organized as follows.
We summarize in Section \ref{sec:galilean}
several properties on Galilean transform which we use,
and then prove Theorem \ref{thm:existence2} in Section \ref{sec:GWP}.
Section \ref{sec:swaves} is devoted to the proof
of Theorems \ref{thm:single} and \ref{thm:multi}.
Then, stability of standing waves is considered in
Section \ref{sec:stable}.

\section{Preliminaries on Galilean transform}\label{sec:galilean}
Let us summarize briefly the properties on the Galilean transform.
\begin{Lem}\label{lem:GT}
Let $\mathcal{G}_\kappa $ be the Galilean transform given in \eqref{eq:GT}.
Then, for any ${\kappa},\,t \in \R$ and ${\bf a},\,{\bf b},\,{\bf a}_1,\,{\bf a}_2,\,{\bf b}_1,\,{\bf b}_2\in\R^d$, we have the following properties.
\begin{enumerate}
\item $\mathcal{G}_{\kappa}(t,0,0)=\mathrm{Id}$ and
$\mathcal{G}_{\kappa}(0,{\bf a},{\bf b})=\mathcal{G}_{0}(0,{\bf a},{\bf b})$.
\item
$	X[\mathcal{G}_{\kappa}(t,{\bf a},{\bf b})u]=X[u]+g_\kappa(t) M[u] $,
$	P[\mathcal{G}_{\kappa}(t,{\bf a},{\bf b})u]=P[u]+g_\kappa^\prime(t) M[u] $.
\item $\mathcal{G}_{\kappa}(t,{\bf a},{\bf b})^{-1}=\mathcal{G}_{\kappa}(t,-{\bf a},-{\bf b})$.
\item $(2i\d_t+\Delta-\ka|x|^2) \mathcal{G}_{\kappa}(t,{\bf a},{\bf b})
	= \mathcal{G}_{\kappa}(t,{\bf a},{\bf b})(2i\d_t+\Delta-\ka|x|^2).$
\item If $N(u)$ is either $F(|u|)u$ or $u\int_{\R^d} K(x-y) |u(y)|^2 dy$,
where $F: \R \to \R$ and $K: \R^d\to \R$, then
$
	N(\mathcal{G}_{\kappa}(t,{\bf a},{\bf b})u)
	=\mathcal{G}_{\kappa}(t,{\bf a},{\bf b}) N(u).
$
\item
$
	\mathcal{U}_\kappa(t) \mathcal{G}_{\kappa}(0,{\bf a},{\bf b})
	= \mathcal{G}_{\kappa}(t,{\bf a},{\bf b})\mathcal{U}_\kappa(t).
$
\item 
	$\mathcal{G}_{\kappa_1}(t,{\bf a}_1,{\bf b}_1)\mathcal{G}_{\kappa_2}(t,{\bf a}_2,{\bf b}_2)u(x)$ is equal to
$$
	e^{\frac{i}2(g_1^\prime\cdot g_2- g_1 \cdot g_2^\prime)(t)}
	e^{-\frac{i}2 (g_1+g_2)(t)\cdot (g_1+g_2)^\prime(t)}
	e^{ix\cdot (g_1+g_2)^\prime(t)}
	u(x-(g_1+g_2)(t)),
$$
where $g_j(t)=g_{\ka_j}(t,{\bf a}_j,{\bf b}_j)$. 
In particular,
$
	\mathcal{G}_{\kappa}(t,{\bf a}_1,{\bf b}_1)\mathcal{G}_{\kappa}(t,{\bf a}_2,{\bf b}_2)
	= e^{i\frac12({\bf a}_1\cdot{\bf b}_2-{\bf a}_2\cdot{\bf b}_1)}\mathcal{G}_{\kappa}(t,{\bf a}_1+{\bf a}_2,{\bf b}_1+{\bf b}_2).
$
\end{enumerate}
\end{Lem}

\begin{proof}
The first two properties are obvious.
The forth and fifth are well known.
The sixth is just a rephrase of Proposition 2.5 in \cite{MHpre}.
One can show the last one by a simple computation.
Note that if $\ka_1=\ka_2$ then
 the Wronskian $g_1^\prime\cdot g_2 - g_1\cdot g_2^\prime$ is 
independent of time and so that
\[
 	\mathcal{G}_{\kappa}(t,{\bf a}_1,{\bf b}_1)\mathcal{G}_{\kappa}(t,{\bf a}_2,{\bf b}_2)
 	= e^{\frac{i}2( {\bf a}_1\cdot {\bf b}_2-{\bf a}_2\cdot {\bf b}_1)} 
 	\mathcal{G}_{\kappa}(t,{\bf a}_1+{\bf a}_2,{\bf b}_1+{\bf b}_2)
\]
follows. Letting ${\bf a}_2=-{\bf a}_1$ and ${\bf b}_2=-{\bf b}_1$,
we obtain the third.
\end{proof}

\begin{Lem}\label{lem:Ssbound}
Let $d\ge1$ and $s\in\R$. There exists a constant $C=C(d,s)$ such that
\[
	\norm{e^{ix\cdot {\bf a}}f}_{\Sigma^s} + \norm{f(\cdot-{\bf a})}_{\Sigma^s}
	\le C \Jbr{|{\bf a}|}^{|s|} \norm{f}_{\Sigma^s}
\]
\end{Lem}
\begin{proof}
This inequality is easily verified when $s=2m$ ($m=0,1,2,\dots$) since
$\norm{f}_{\Sigma^s}$ is written as $\Lebn{(-\frac12 \Delta + \frac12 |x|^2)^mf}2$.
Then, the case $s\ge 0$ immediately follows by interpolation.
For $s<0$, we use a duality argument.
\end{proof}

\begin{Lem}\label{lem:Sslimit}
Let $s\in\R$ and $f \in \Sigma^s$. Then,
$\norm{e^{ix\cdot {\bf a}}f-f}_{\Sigma^s} + \norm{f(\cdot-{\bf a})-f}_{\Sigma^s}
\to 0$ as $|{\bf a}|\to0$.
\end{Lem}
\begin{proof}
Set $f_N = \sum_{|\underline{n}|\le N} (f,\Omega_{\underline{n}})\Omega_{\underline{n}}$.
Then, $f_N \to f$ in $\Sigma^s$ as $N\to\I$.
By Lemma \ref{lem:Ssbound},
for any $\eps>0$ there exists $N_0$ such that $\norm{f-f_{N_0}}_{\Sigma^s}\le
\eps/3$ and 
\[
	\norm{e^{ix\cdot {\bf a}}(f-f_{N_0})}_{\Sigma^s} + \norm{(f-f_{N_0})(\cdot-{\bf a})}_{\Sigma^s}
	\le \frac{\eps}3
\]
for any $|{\bf a}|\le 1$.
Take $m\in \N$ so that $s\le 2m$.
For any fixed multi-index $\underline{n}$,
we deduce from Lebesgue's convergence theorem that
\begin{multline*}
	\Lebn{\(-\frac12 \Delta + \frac12 |x|^2\)^m(e^{ix\cdot {\bf a}}\Omega_{\underline{n}}-\Omega_{\underline{n}})}2 \\
	+ \Lebn{\(-\frac12 \Delta + \frac12 |x|^2\)^m(\Omega_{\underline{n}}(\cdot-{\bf a})-\Omega_{\underline{n}})}2
\to 0
\end{multline*}
as $|{\bf a}|\to0$.
Since
$\norm{u}_{\Sigma^s}\le \norm{u}_{\Sigma^{2m}} =\Lebn{(-\frac12 \Delta + \frac12 |x|^2)^{m}u}2$ for $u\in \Sigma^s$
and since $f_{N_0}$ is a finite combination of $\Omega_{\underline{n}}$,
there exists $ \delta \in(0, 1]$ such that if $|{\bf a}|\le \delta$ then
\[
	\norm{e^{ix\cdot {\bf a}}f_{N_0}-f_{N_0}}_{\Sigma^s} + \norm{f_{N_0}(\cdot-{\bf a})-f_{N_0}}_{\Sigma^s}
	\le \frac{\eps}3.
\]
\end{proof}
From Lemma \ref{lem:Ssbound} and \ref{lem:Sslimit},
we obtain the following.
\begin{Lem}\label{lem:GTcont}
Let $s\in\R$. For any $t,\ka\in\R$ and ${\bf a},{\bf b}$,
the Galilean transform $\mathcal{G}_\ka(t,{\bf a},{\bf b})$ is a bounded linear map
from $\Sigma^s$ to itself. The operator norm has the following estimate:
\[
	\norm{\mathcal{G}_\ka(t,{\bf a},{\bf b})}_{\mathcal{B}(\Sigma^s)}\le
	\begin{cases}
	(1+C (|{\bf a}|+|{\bf b}|) )^{|s|}& \ka>0,\\
	(1+C(|{\bf a}|t + |{\bf a}|+|{\bf b}| ) )^{|s|} & \ka=0, \\
	(1+C(|{\bf a}|+|{\bf b}|)e^{|t\sqrt{\ka}|})^{|s|} & \ka<0.
	\end{cases}
\]
In particular, for each $u\in \Sigma^s$, a mapping $\R^{1+1+d+d} \ni(t,\ka,{\bf a},{\bf b}) \mapsto 
\mathcal{G}_\ka(t,{\bf a},{\bf b})u \in \Sigma^s$ is continuous.
\end{Lem}

\section{Global well-posedness of \eqref{eq:nh2}}\label{sec:GWP}
In this section, we prove Theorem \ref{thm:existence2}.

\begin{Lem}\label{lem:XP}
Let $u\in\Sigma^{1/2}$.
Then, $X[u]$ and $P[u]$ are well-defined.
Further, if $\{u_n\}_n \subset \Sigma^{1/2}$ converges
to $u$ in $\Sigma^{1/2}$ as $n\to\I$ then
$X[u_n]$ and $P[u_n]$ converges to $X[u]$ and $P[u]$, respectively, as $n\to\I$.
\end{Lem}
\begin{proof}
For each $i\in[1,d]$,
one verifies that
\begin{equation}\label{eq:XPformula}
\begin{aligned}
	\int_{\R^d} x_i |u(x)|^2 dx = {}&
	\mathrm{Re}\sum_{\underline{n}} \sqrt{2(n_i+1)}a_{\underline{n}}
	\overline{a_{\underline{n}+{\bf e}_i}}, \\
	\mathrm{Im}\int_{\R^d} \overline{u(x)}\d_{i}u(x) dx ={}& 
	\mathrm{Im}\sum_{\underline{n}} \sqrt{2(n_i+1)}a_{\underline{n}}
	\overline{a_{\underline{n}+{\bf e}_i}},
\end{aligned}
\end{equation}
where $a_{\underline{n}} = (u,\Omega_{\underline{n}})$ and
${\bf e}_i$ is a multi-index such that the $i$-th component is
one and the others are zero. Therefore,
$|X[u]|+|P[u]| \le C \norm{u}_{\Sigma^{1/2}}$ holds.
Convergence part follows from
\[	
	|X[u]-X[v]| + |P[u]-P[v]|
	\le C(\norm{u}_{\Sigma^{1/2}}+\norm{v}_{\Sigma^{1/2}})\norm{u-v}_{\Sigma^{1/2}}.
\]
\end{proof}

\begin{Lem}\label{lem:L2}
Let $V(t):\R\to\R$ be a continuous function of time.
If $u\in C(\R;L^2)\cap C^1(\R;\Sigma^{-2})$ solve
\[
	2i u_t +\Delta u - V(t)|x|^2u =0
\]
in $\Sigma^{-2}$ sense. then $\Lebn{u(t)}2$ is constant. 
\end{Lem}
\begin{proof}
By an abstract theory by Kato \cite{TKato70},
there exists a unique propagator $\{U(t,s)\}_{t,s \in \R}$ with
the following properties;
\begin{enumerate}
\item  $U(t,s)$ is unitary on $\Sigma^{-2}$ and strongly continuous in $s,t$;
\item  $U(t,t)=1$ and $U(t,s)=U(t,r)U(r,s)$ for any $t,s,r\in\R$;
\item  
$U(t,s){L^2}\subset L^2$.
Further, $\left.U(t,s)\right|_{L^2}$ is unitary on $L^2$ and strongly continuous in $s,t$;
\item  $\frac{d}{ds}U(t,s)y=U(t,s)i(-\frac12\Delta+\frac{V(s)}2|x|^2)y \in \Sigma^{-2}$
for any $y\in L^2$;
\item  $\frac{d}{dt}U(t,s)y=i(\frac12\Delta-\frac{V(t)}2|x|^2)U(t,s)y \in \Sigma^{-2}$
for any $y\in L^2$
\end{enumerate}
(See Theorems 4.1 and 5.1 and Remarks 5.3 and 5.4 of \cite{TKato70}).
Set $w(s)=U(t,s)u(s)$. By the forth property of $\{U(t,s)\}$
and by assumption on $u$,
we obtain
\[
	\frac{d}{ds}w(s) = U(t,s)
	i\(-\frac12\Delta+\frac{V(s)}2|x|^2\) u(s)
	+U(t,s) u_t(s) =0 \IN \Sigma^{-2} 
\]
for any $0<s<t$.
Hence, $u(t)=w(t)=w(0)=U(t,0)u(0)$ for $t>0$.
This holds also for $t<0$ by the same argument.
Thus, the result is true because of the third property of $\{U(t,s)\}$.
\end{proof}
We are now in a position to prove Theorem \ref{thm:existence2}.
\begin{proof}[Proof of Theorem \ref{thm:existence2}]
Let $v_0 \in \Sigma^{1/2}$ and let $M$, ${\bf a}$, ${\bf b}$, and $\ka$ be as in
the statement of Theorem \ref{thm:existence2}.
We shall first prove that $v(t)$ given by the first line of
\eqref{eq:mf2} is a solution
of \eqref{eq:nh2}.
Set $w(t)=\mathcal{U}_{\ka}(t)\mathcal{G}_\ka(0,{\bf a},{\bf b})^{-1} v_0$.
It holds that  $w$ is in $ C(\R;\Sigma^{1/2})\cap C^1(\R;\Sigma^{-3/2})$ and
solves $2i w_t =- \Delta w + \ka |x|^2 w$ in $\Sigma^{-3/2}$ sense
with $w(0)=\mathcal{G}_\ka(0,{\bf a},{\bf b})^{-1} v_0$.
Moreover, $M[w(t)]=M$ and $X[w(t)]\equiv P[w(t)] \equiv 0$.
Now, we set $z(t)=\mathcal{G}_\la (t,{\bf a},{\bf b}) w(t)$.
Then, $z(t)\in C(\R;\Sigma^{1/2})$ by Lemma \ref{lem:GTcont}.
Similarly, $z(t)\in C^1(\R;\Sigma^{-3/2})$ follows from
\[
	2i z_t = (|g^\prime_\la(t)|^2+g_\la(t)\cdot g_\la^{\prime\prime}(t))z
	-2 g^{\prime\prime}_\la(t)\cdot x z - 2i g^\prime_\la(t)\cdot \mathcal{G}_\la(t)\nabla w + \mathcal{G}_\la(t) 2i w_t.
\]
We have $M[z(t)]=M$.
It also holds from the second property of Lemma \ref{lem:GT}
that $X[z(t)]=M g_\la(t)$. 
The forth property of Lemma \ref{lem:GT} implies
$2i z_t + \Delta z - \la|x|^2 z =\mathcal{G}_\la (t,{\bf a},{\bf b}) (2i w_t + \Delta w - \la|x|^2 w)=\mathcal{G}_\la (t,{\bf a},{\bf b})\eta M |x|^2 w $.
Since
\begin{align*}
	\mathcal{G}_\la(t) \eta M |x|^2 w={}& \eta M|x|^2z
	- 2\eta M g_\la(t)\cdot x z + \eta M|g_\la(t)|^2 z \\
	 ={}& \eta |x|^2M[z(t)] z - 2\eta x \cdot X[z(t)] z + 2\Phi^\prime(t) z \\
	 ={}& \eta \int(|x-y|^2-|y|^2)|z(t,y)|^2 dy z + 2\Phi^\prime(t) z,
\end{align*}
we conclude that $v(t):=z(t)e^{-i\Phi(t)}$ is a solution of \eqref{eq:nh2}.
Statements on energy conservation follows from Lemma \ref{lem:energy}, below.

The data-to-solution map $v_0\mapsto v(t)$ is continuous
as a map from $\Sigma^{1/2}$ to $L^\I_{\mathrm{loc}}(\R; \Sigma^{1/2})$
due to Lemmas \ref{lem:GTcont} and \ref{lem:XP}.

Let us proceed to the uniqueness.
Assume $\widetilde{v}(t) \in C(\R;\Sigma^{1/2})\cap C^1(\R;\Sigma^{-3/2})$ solves
\eqref{eq:nh2} in $\Sigma^{-3/2}$ sense with the same initial condition
$\widetilde{v}(0)=v_0$.
Since 
$M[\widetilde{v}(t)]$ and $X[\widetilde{v}(t)]$ are continuous,
there exists an $\R^d$-valued function $H(t)\in C^2(\R)$ such that
\begin{equation}\label{eq:pf2-2}
H^{\prime\prime}(t)=-(\la +\eta M[\widetilde{v}(t)])H(t)+\eta X[\widetilde{v}(t)]
\end{equation}
and $H(0)=H^\prime(0)=0$.
Define $\widetilde{w}(t)$ by
\begin{equation}\label{eq:pf2-3}
	\widetilde{v}(t,x)
	= e^{-i\widetilde{\Phi}(t)}e^{-\frac{i}2 H(t)\cdot H^\prime(t)}e^{ix\cdot H^\prime(t)}\widetilde{w}(t,x-H(t)),
\end{equation}
where
$
	\widetilde{\Phi}(t) = -\frac{\eta}2 \int_0^t H(s)\cdot X[\widetilde{v}(s)] ds.
$
It is easy to see $\widetilde{w}\in C(\R;\Sigma^{1/2})\cap C^1(\R;\Sigma^{-3/2})$ since $H(t)\in C^2(\R)$.
We have $M[\widetilde{w}]=M[\widetilde{v}]$ and
$X[\widetilde{w}]= X[\widetilde{v}]-M[\widetilde{v}] H$.
%
%
%
%
By a computation with an identity
\begin{align*}
	&\eta \int (|x-y|^2-|y|^2) |\widetilde{v}(y)|^2 dy \widetilde{v} \\={}&
	\eta |x|^2 M[\widetilde{v}]\widetilde{v}
	-2\eta x\cdot X[\widetilde{v}] \widetilde{v}\\
	={}&
	e^{-i\widetilde{\Phi}(t)}e^{-\frac{i}2 H\cdot H^\prime}e^{ix\cdot H^\prime}
	(\eta M[\widetilde{v}] |x|^2 \widetilde{w})(t,x-H(t)) 
	+ 2\eta M[\widetilde{v}] H\cdot x\widetilde{v}
	- \eta M[\widetilde{v}] |H|^2 \widetilde{v}
	\\&{}-2\eta  X[\widetilde{v}] \cdot x \widetilde{v}
\end{align*}
and \eqref{eq:pf2-2}, one verifies that $\widetilde{w}$ solves
$
	2i \widetilde{w}_t + \Delta \widetilde{w} - (\la + \eta M[\widetilde{w}]) |x|^2 \widetilde{w} =0
$
in $\Sigma^{-3/2}$ with $\widetilde{w}(0)=v_0$.
Since $V(t):=\la + \eta M[\widetilde{w}(t)]$ is continuous in time by assumption,
$M[\widetilde{w}(t)]\equiv M[v_0]$
follows from Lemma \ref{lem:L2}.
Set $M:=M[v_0]$ and $\ka=\la + \eta M$.
$\widetilde{w}$ is written as $\widetilde{w}(t)=\mathcal{U}_\ka(t)v_0$ and
so $X[\widetilde{w}(t)]=M g_\ka(t,{\bf a},{\bf b})$.
The equation which $H$ solves now
becomes $H^{\prime\prime}(t)= -\la H(t) + \eta M g_\ka(t) $ and $H(0)=H^\prime(0)=0$.
Solving this, we see
$H(t)=g_\la(t)-g_\ka(t)$.
Consequently, $X[\widetilde{v}(t)]=M g_\la(t)$ and
\[
	\widetilde{\Phi}(t) = -\frac{\eta M}2 \int_0^t (g_\la-g_\ka)(s)\cdot
	g_\la(s) ds
	=\Phi(t) + \frac{1}2(g_\la (t)\cdot g_\ka^\prime(t)-g_\la^\prime(t) \cdot g_\ka(t)).
\]
Substituting these formulas to \eqref{eq:pf2-3} and using 
the third and the seventh properties of Lemma \ref{lem:GT},
we see that $\widetilde{v}(t)$ coincides with $v(t)$.
\end{proof}

\section{Standing waves}\label{sec:swaves}
We now turn to the proof of Theorems \ref{thm:single}
and \ref{thm:multi}.
\begin{Lem}\label{lem:Psi}
Let $u_0 \in \Sigma$ and
$\ka=\la+\eta M>0$. The phase $\Psi(t)$ in Theorem \ref{thm:existence} is written as
\[
	\Psi(t) = \frac{\eta}{2\kappa^{\frac12}} \(\Lebn{A_\ka w_0}2^2 
	+ \Lebn{A_\ka^\dagger w_0}2^2\) 
	t 
	+\frac{\eta}{ \ka} \sin (\sqrt\ka t) \Re \(  e^{-i\sqrt\ka t} (A_\ka w_0, A_\ka^\dagger w_0) \),
\] 
where $w_0=\mathcal{G}_\ka(0,{\bf a}, {\bf b})^{-1}u_0$,
$A_\kappa= \frac1{\sqrt2} (\ka^{1/4}x + \ka^{-1/4}\nabla )$,
and $A_\kappa^\dagger=\frac1{\sqrt2}( \ka^{1/4}x - \ka^{-1/4}\nabla )$.
\end{Lem}
\begin{proof}
It is well known that 
\[
	x \mathcal{U}_\ka(t)= \mathcal{U}_\ka(t) 
	\(\cos(\sqrt\ka t) x -i\frac{\sin (\sqrt\ka t)}{\sqrt\ka}\nabla \).
\]
By definitions of $A_\ka$ and $A_\ka^\dagger$, the right hand side is written as
\[
	\frac1{\sqrt{2\ka^{1/2}}} (e^{-i\sqrt\ka t}A_\ka + e^{i\sqrt\ka t}A_\ka^\dagger).
\]
It therefore holds that
\begin{align*}
	\Lebn{x\mathcal{U}_\ka(t)w_0}2^2
	={}& \frac1{2\ka^{1/2}}\(\Lebn{A_\ka w_0}2^2 + \Lebn{A_\ka^\dagger w_0}2^2\) \\
	&{}+ \frac1{\ka^{1/2}} \Re \(e^{-2i\sqrt\ka t}(A_\ka w_0, A_\ka^\dagger w_0)\). 
\end{align*}
Thus, integration in time gives us the desired result.
\end{proof}
\begin{proof}[Proof of Theorem \ref{thm:single}]
Let  $u_0=M^\frac12 \mathcal{G}_\ka(0,{\bf a}_1, {\bf b}_1)
 \Omega_{\underline{n},\ka}(x)$.
Then,
$M[u_0]=M$ holds.
$P[u_0]=M{\bf a}_1$ and $X[u_0]=M{\bf b}_1$ follow
 from the second property of Lemma \ref{lem:GT}
since $P[\Omega_{\underline{n},\ka}]=X[\Omega_{\underline{n},\ka}]=0$.
One sees that $\mathcal{G}_\ka(0,{\bf a}_1,{\bf b}_1)^{-1}u_0=\Omega_{\underline{n},\ka}$.
Therefore, the formula \eqref{eq:mf} tells us that
\[
	u(t) = e^{-i\Psi(t)}\mathcal{G}_\la(t,{\bf a}_1,{\bf b}_1)
	\mathcal{U}_\ka(t) \Omega_{\underline{n},\ka}
	= e^{-i\Psi(t)}\mathcal{G}_\la(t,{\bf a}_1,{\bf b}_1)
	e^{-\frac{i}2\ka^{\frac12}(|\underline{n}|+\frac{d}2)t} \Omega_{\underline{n},\ka}.
\]
Moreover, we deduce from Lemma \ref{lem:Psi} that
$\Psi(t)=\frac{\eta M}{2\ka^{1/2}}(|\underline{n}|+\frac{d}2)t$.
The second part is shown by a similar argument.
We only note that
$\Phi(t)$ is proportional to $t$ if and only if $|{\bf a}_1|^2=\la |{\bf b}_1|^2$,
and that $\Phi(t)=-\frac{\eta M |{\bf b_1}|^2}2 t$ under this condition.
\end{proof}
\begin{proof}[Proof of Theorem \ref{thm:multi}]
Let ${\bf a}$ and ${\bf b}$ be as in assumption.
We only consider \eqref{eq:nh}.
Since the Galilean transform is a linear map,
\[
	\mathcal{G}_\ka(0,{\bf a},{\bf b})^{-1}u_0
	= \mu \sum_{j=1}^L \widetilde{\alpha}_j \mathcal{G}_\ka(0,{\bf a}_j-{\bf a}, {\bf b}_j-{\bf b})
	\Omega_{\underline{n}_j,\ka},
\]
where we have applied the third and the seventh properties of
Lemma \ref{lem:GT}. Thanks to the sixth properties of Lemma \ref{lem:GT},
we see
\[
	\mathcal{U}_\ka(t)\mathcal{G}_\ka(0,{\bf a},{\bf b})^{-1}u_0
	= \mu \sum_{j=1}^L \widetilde{\alpha}_j \mathcal{G}_\ka(t,{\bf a}_j-{\bf a}, {\bf b}_j-{\bf b})
	(e^{-\frac{i}2\ka^{\frac12}(|\underline{n}_j|+\frac{d}2)t} \Omega_{\underline{n}_j,\ka}).
\]
Now the result is obvious by \eqref{eq:mf}.
\end{proof}
\section{Stable excited states}\label{sec:stable}

This section is devoted to the proof of Theorems \ref{stability} and \ref{thm:main}.
We prepare several lemmas for the proof of theorems.

\begin{Lem}\label{lem:ab}
Let $\{u_n\}\subset \Sigma^{1/2}$.
If $||u_n-M^{1/2}\Om_{\underline{n},\ka}||_{\Sigma^{1/2}}\rightarrow 0$, then $P[u_n]\rightarrow 0$ and $X[u_n]\rightarrow0$.
\end{Lem}

\begin{proof}
This immediately follows from Lemma \ref{lem:XP} and the fact that 
$P[M^{1/2}\Om_{\underline{n},\ka}]=0$ and $X[M^{1/2}\Om_{\underline{n},\ka}]=0$.
\end{proof}

\begin{Lem}\label{lem:g}
Let $\al\geq 0$, then $\sup_{t>0}|g_{\al}'(t,{\bf a}, {\bf b})|\rightarrow 0$ as $|{\bf a}|+|{\bf b}|\rightarrow 0$.
\end{Lem}

\begin{proof}
For the case $\al=0$, we have $g_0(t,{\bf a},{\bf b})={\bf a}t+{\bf b}$ and for the case $\al>0$, we have $g_{\al}(t)=a^{-1/2}{\bf a}\sin \al^{1/2}t+{\bf b}\cos \al^{1/2}t$.
Thus, the conclusion of the lemma is obvious.
\end{proof}

We prove Theorem \ref{stability} by using Lemmas \ref{lem:Ssbound}, \ref{lem:Sslimit}, \ref{lem:ab}, and \ref{lem:g}.
\begin{proof}[Proof of Theorem \ref{stability}]
We consider standing waves of \eqref{eq:nh}.
Fix $\underline{n}$ and set $Q_{\underline{n},M}=e^{-i\om_1t/2}M^{1/2}\Om_{\underline{n},\ka}$.
Let $s\geq 1$ and take $u_0\in\Sigma^s$ with $\Lebn{u_0}2^2=M$.
Since $\{\Om_{\underline{n},\ka}\}$ is a CONS of $L^2$,
$u_0$ is  decomposed as
\[
	u_0=\sum_{\underline{m}}a_{\underline{m}}\Om_{\underline{m},\ka}.
\]
By using the formula (\ref{eq:mf}), 
one sees that
\begin{align*}
	u(t)={}& e^{i\tilde{\Psi}(t)}e^{ix\cdot (g_{\ka}'(t,-{\bf a},-{\bf b})+g_{\la}'(t,{\bf a},{\bf b}))}\\
	&{} \times\sum_{\underline{m}}e^{-it\ka^{1/2}\(|\underline{m}|+\frac{d}{2}\)}a_{\underline{m}}\Om_{\underline{m},\ka}(x-g_{\ka}(t,-{\bf a},-{\bf b})-g_{\la}(t,{\bf a},{\bf b})),
\end{align*}
where $\tilde{\Psi}(t)$ is a function independent of $x$.
Therefore, we have
\begin{equation}\label{eq:u}
\begin{aligned}
	&e^{i\tilde{\Phi}(t)}u(t,x+g_{\ka}(t,-{\bf a},-{\bf b})+g_{\la}(t,{\bf a},{\bf b}))\\
	&{}=e^{ix\cdot (g_{\ka}'(t,-{\bf a},-{\bf b})+g_{\la}'(t,{\bf a},{\bf b}))}a_{\underline{n}}\Om_{\underline{n},\ka}\\
	&\qquad +e^{ix\cdot (g_{\ka}'(t,-{\bf a},-{\bf b})+g_{\la}'(t,{\bf a},{\bf b}))}\sum_{\underline{m}\neq \underline{n}}e^{-it\ka^{1/2}\(|\underline{m}|-|\underline{n}|+\frac{d}{2}\)}a_{\underline{m}}\Om_{\underline{m},\ka},
\end{aligned}
\end{equation}
 where $\tilde{\Phi}$ is another function independent of $x$.
The difference between $Q_{\underline{n},M}$ and (\ref{eq:u}) 
is hence calculated as follows;
\begin{equation}\label{eq:estdiff}
\begin{aligned}
	&\norm{Q_{\underline{n},M}-e^{i\tilde{\Phi}(t)+i\om_1 t/2}u(t,x+g_{\ka}(t,-{\bf a},-{\bf b})+g_{\la}(t,{\bf a},{\bf b}))}_{\Sigma^s}\\
	&{}\leq |M^{1/2}-a_{\underline{n}}|\norm{\Om_{\underline{n},\ka}}_{\Sigma^s}
	+|a_{\underline{n}}|\norm{(1-e^{ix\cdot (g_{\ka}'(t,-{\bf a},-{\bf b})+g_{\la}'(t,{\bf a},{\bf b}))})\Om_{\underline{n},\ka}}_{\Sigma^s}\\
	&\qquad +\norm{e^{ix\cdot (g_{\ka}'(t,-{\bf a},-{\bf b})+g_{\la}'(t,{\bf a},{\bf b}))}\sum_{\underline{m}\neq \underline{n}}a_{\underline{m}}\Om_{\underline{m}}}_{\Sigma^s}.
\end{aligned}
\end{equation}
It is obvious from Lemma \ref{lem:Ssbound}
that the first term and the third term in \eqref{eq:estdiff} are small when $\norm{u_0-M^{1/2}\Om_{\underline{m},\ka}}_{\Sigma^s}$ is sufficiently small.
We deduce from Lemmas \ref{lem:ab} and \ref{lem:g} that
if $\norm{u_0-M^{1/2}\Om_{\underline{m},\ka}}_{\Sigma^s}$ is small then
$\sup_{t>0}|g_{\ka}'(t,-{\bf a},-{\bf b})+g_{\la}'(t,{\bf a},{\bf b}))|\ll 1$.
It then follows from Lemma \ref{lem:Sslimit} that the second term in \eqref{eq:estdiff}
is also small.

For the case $\Lebn{u_0}2^2\neq M$, we can show that $u(t)$ is near $Q_{\underline{n},M'}$ where $M'=\Lebn{u_0}2^2$.
Since for $u_0$ sufficiently near $Q_{\underline{n},M}$, we have $|M-M'|\ll 1$.
So, $Q_{\underline{n},M}$ and $Q_{\underline{n},M'}$ is also near to each other
(up to phase).
Therefore, we have the conclusion.
The proof for the standing waves of \eqref{eq:nh2} is similar.
We use \eqref{eq:mf2} instead of \eqref{eq:mf}.
\end{proof}

Next lemma is concerned with energy.
By this lemma, one sees that, in Theorem \ref{stability},
the standing wave solution is a ground state if $\underline{n}=0$
and an excited state otherwise.
This is an immediate consequence of the representation of the solutions,
so we omit details of the proof.
\begin{Lem}\label{lem:energy}
Let $d\ge1$, $\la,\eta\in\R$, and $u_0,v_0\in \Sigma$ with $M[u_0]=M[v_0]=M$.
Let $\ka=\la+\eta M$.
Let $u(t)$ and $v(t)$ be solutions to \eqref{eq:nh} and \eqref{eq:nh2}
given in Theorems \ref{thm:existence} and \ref{thm:existence2}, respectively.
Then, the energy $E[u(t)]$ (resp. $E[v(t)]$) is conserved and given by
\[
	\(\(-\frac12\Delta + \frac{\ka}2|x|^2\)w_0,w_0\)
	+ \frac{M}2(|{\bf a}|^2+\la |{\bf b}|^2),
\]
where $w_0=\mathcal{G}_\ka(0,{\bf a},{\bf b})^{-1}u_0$
(resp.\ $w_0=\mathcal{G}_\ka(0,{\bf a},{\bf b})^{-1}v_0$).
Further, suppose $\ka>0$.
Set $$e(M):=\inf_{u\in\Sigma, M[u]=M} E[u].$$
The following holds.
\begin{enumerate}
\item If $\l>0$ then $E[u_0]=e(M)$ holds if and only if
$w_0=M^{1/2}\Om_{0,\ka}$ and ${\bf a}={\bf b}=0$.
\item If $\l=0$ then $E[u_0]=e(M)$ holds if and only if
$w_0=M^{1/2}\Om_{0,\ka}$ and ${\bf a}=0$.
\item If $\l<0$ then $e(M)=-\I$.
\end{enumerate}
\end{Lem}

We now move to the proof of Theorem \ref{thm:main}.
Before the proof, let us introduce some notations.
In what follows, we shall consider  \eqref{eq:nh} in the following cases.
\begin{enumerate}
\item[$(\mathrm{I})$]
$d=1$, $\la=0$ and $\eta=1$,
\item[$(\mathrm{II})$]
$d=1$, $\la=2$ and $\eta=-1$.
\end{enumerate}
Denote the equation \eqref{eq:nh} for the cases $\mathrm{I}$ and $\mathrm{II}$ as $\eqref{eq:nh}_\mathrm{I}$ and $\eqref{eq:nh}_{\mathrm{II}}$, respectively.
By Theorem \ref{thm:single}, the following functions are solutions of $\eqref{eq:nh}_\mathrm{I}$ and $\eqref{eq:nh}_{\mathrm{II}}$.
\begin{eqnarray*}
Q_{\mathrm{I},n}(t,x;M)&:=&\mathrm{exp}\(-\frac{3}{2}M^{1/2}\(n+\frac{1}{2}\)t\)M^{\frac{1}{2}+\frac{1}{8}}\Om_n(M^{1/4}x),\\
Q_{\mathrm{II},n}(t,x;M)&:=&\mathrm{exp}\(-\(n+\frac{1}{2}\)\(\ka^{1/2}-\frac{1}{2}M\ka^{-1/2}\)t\)M^{\frac{1}{2}}\ka^{\frac{1}{8}}\Om_n(\ka^{1/4}x),
\end{eqnarray*}
where $\ka=2-M$.
We can also write
\begin{eqnarray*}
Q_{\mathrm{I},n}=e^{-\frac{\om_{\mathrm{I},n}(M)}{2}t}M^{\frac{1}{2}+\frac{1}{8}}\Om_n(M^{1/4}x),\\
Q_{\mathrm{II},n}=e^{-\frac{\om_{\mathrm{II},n}(M)}{2}t}M^{\frac{1}{2}}\ka^{\frac{1}{8}}\Om_n(\ka^{1/4}x),
\end{eqnarray*}
where
\begin{eqnarray*}
\om_{\mathrm{I},n}(M)&:=&3M^{1/2}\(n+\frac{1}{2}\),\\
\om_{\mathrm{II},n}(M)&:=&2\(n+\frac{1}{2}\)\(\ka^{1/2}-\frac{1}{2}M\ka^{-1/2}\).
\end{eqnarray*}
Let $S_{\om,\mathrm{I}}$, $S_{\om,\mathrm{II}}$ be the action of $\eqref{eq:nh}_\mathrm{I}$ and $\eqref{eq:nh}_{\mathrm{II}}$, respectively.
Further, set
\begin{eqnarray*}
d_\mathrm{I}(\om)&:=&S_{\om,\mathrm{I}}\(M^{\frac{1}{2}+\frac{1}{8}}\Om_n(M^{1/4}x)\),\\
d_\mathrm{II}(\om)&:=&S_{\om,\mathrm{II}}\(M^{\frac{1}{2}}\ka^{\frac{1}{8}}\Om_n(\ka^{1/4}x)\).
\end{eqnarray*}
Since $Q_{\mathrm{I},n}$ and $Q_{\mathrm{II},n}$ are stable (Theorem \ref{stability}),
it suffices to show the following lemma.

\begin{Lem}\label{lem:specq}
Let $n$ be an even integer.
Then, we have the following.
\begin{enumerate}
\item[$(\mathrm{i})$]
$\left.d_{\mathrm{I}}''(\om)\right|_{\om=\om(1)}<0$ and $n\(\left.S_{\om,\mathrm{I}}''(\Om_n)\right|_{L^2_{r}}\)=n$.
\item[$(\mathrm{ii})$]
$\left.d_{\mathrm{II}}''(\om)\right|_{\om=\om(1)}>0$ and $n\(\left.S_{\om,\mathrm{II}}''(\Om_n)\right|_{L^2_{r}}\)=n+1$.
\end{enumerate}
\end{Lem}







\begin{proof}[Proof of Lemma \ref{lem:specq}]
We first calculate $d''$ for the cases (I) and (II).
Since $d(\om)=S_{\om}(\phi_{\om})$, one obtains $d'(\om)=-Q(\phi_{\om})$.
Thus, $d''(\om)=-\frac{1}{2}\frac{dM}{d\om}$.
For the case $(\mathrm{I})$, we have $\om_{\mathrm{I},n}(M)=3M^{1/2}\(n+\frac{1}{2}\)$.
Differentiate $\om_{\mathrm{I},n}$ with respect to $M$ to yield
$$\left.\frac{d\om_{\mathrm{I},n}}{dM}\right|_{M=1}=\frac{3}{2}\(n+\frac{1}{2}\)>0,$$
which implies $\left.d_{\mathrm{I}}''(\om)\right|_{M=1}<0$.
%
On the other hand, we have $\om_{\mathrm{II},n}(M)=\(2n+1\)(2-M)^{-1/2}\(2-\frac{3}{2}M\)$
since $\ka=2-M$.
One hence verifies that
$$\left.\frac{d\om_{\mathrm{II},n}}{dM}\right|_{M=1}=-\frac{5}{2}\(n+\frac{1}{2}\)<0,$$
which yields $\left.d_{\mathrm{II}}''(\om)\right|_{M=1}>0$.
 
We next calculate the spectrum of $S_{\om}''(\Om_n)$ for the cases (I) and (II).
We express $S_{\om}''(\Om_n)$ by a $2\times 2$ matrix $(L_{i,j})_{i,j=1,2}$ as follows
\begin{eqnarray*}
\(
\begin{matrix}
\mathrm{Re}\(S_{{\mathrm{I}},\om}''(\Om_n)u\) \cr
\mathrm{Im}\(S_{{\mathrm{I}},\om}''(\Om_n)u\) \cr
\end{matrix}
\)
=\left(
  \begin{matrix}
        L_{11} & L_{12} \cr
        L_{21} & L_{22} \cr
  \end{matrix}
\right)
\left(
  \begin{matrix}
        \mathrm{Re}\,u \cr
        \mathrm{Im}\,u \cr
  \end{matrix}
\right).
\end{eqnarray*}
For the case $(\mathrm{I})$, we have $E=\frac{1}{2}\int_{\R}|\nabla u|^2+\frac{1}{4}\int_{\R}\int_{\R}|x-y|^2|u(x)|^2|u(y)|^2$.
So, we have
\begin{eqnarray*}
&&S_{{\mathrm{I}},\om}''(\Om_n)\\
&&=\left(
  \begin{matrix}
        -\De+|x|^2+2(|x|^2*(\Om_n\cdot))\Om_n-2\(n+\frac{1}{2}\) & 0 \cr
        0 & -\De+|x|^2-2\(n+\frac{1}{2}\) \cr
  \end{matrix}
\right).
\end{eqnarray*}
It is easy to see that
$L_{22}=-\De+|x|^2-2\(n+\frac{1}{2}\)$ has $n/2$ negative eigenvalues and one zero eigenvalue.
We hence investigate the spectrum of $L_{11}=-\De+|x|^2+2(|x|^2*(\Om_n\cdot))\Om_n-2\(n+\frac{1}{2}\)=L_{22}+2(|x|^2*(\Om_n\cdot))\Om_n$.
Let $h=\sum_{n\in 2\N\cup\{0\}}a_l\Om_l$.
Then, it holds that
\begin{align*}
 &\Jbr{\(|x|^2*(\Om_n h)\)\Om_n, h} \\
&{}= \int_{\R}\int_{\R}|x-y|^2\Om_n(x)\Om_n(y)h(x)h(y)\,dydx\\
&{}= 2\int_{\R}|x|^2\Om_n(x)h(x)\,dx\int_{\R}\Om_n(y)h(y)\,dy\\
&{}= \frac{\sqrt{(n+1)(n+2)}}{2}a_{n+2}a_n+\(n+\frac{1}{2}\)a_n^2+\frac{\sqrt{n(n-1)}}{2}a_na_{n-2}
\end{align*}
Therefore, the term $2(|x|^2*(\Om_n\cdot))\Om_n$ affects only the frame $\{\Om_{n-2}, \Om_n, \Om_{n+2}\}$.
In this frame, the action of $L_{11}$ can be
represented as the following matrix.
\begin{eqnarray*}
A_{\mathrm{I}}=(A_{ij})_{i,j=1,2,3}=
\left(
  \begin{matrix}
        -2 & \frac{\sqrt{n(n-1)}}{4} & 0\cr
        \frac{\sqrt{n(n-1)}}{4} & n+\frac{1}{2} & \frac{\sqrt{(n+1)(n+2)}}{4} \cr
        0 & \frac{\sqrt{(n+1)(n+2)}}{4} & 2
  \end{matrix}
\right),
\end{eqnarray*}
where $L_{11}(a_1\Om_{n-2}+a_2\Om_{n}+a_2\Om_{n+2})=\(\sum_{j=1,2,3}A_{ij}a_j\)\Om_{n-4+2i}$.
The characteristic polynomial of $A_{\mathrm{I}}$ becomes
\begin{eqnarray*}
F_{\mathrm{I}}(\la)=\la^3-\(n+\frac{1}{2}\)\la^2-\frac{1}{8}\(n^2+n+33\)\la+\frac{7}{2}\(n+\frac{1}{2}\).
\end{eqnarray*}
Since $F_{\mathrm{I}}(0)>0>F_{\mathrm{I}}(n)$ holds for any even integer $n$,
 $A_{\mathrm{I}}$ has one negative eigenvalue and two positive eigenvalues.
Therefore, we have the conclusion.

Now we consider the case $(\mathrm{II})$.
As in the case (${\mathrm{I}}$), since $E=\frac{1}{2}\int_{\R}|\nabla u|^2+\int_{\R}|x|^2|u|^2-\frac{1}{4}\int_{\R}\int_{\R}|x-y|^2|u(x)|^2|u(y)|^2$,
one obtains
\begin{eqnarray*}
&&S_{{\mathrm{II}},\om}''(\Om_n)\\
&&=\left(
  \begin{matrix}
        -\De+|x|^2+2(|x|^2*(\Om_n\cdot))\Om_n-2\(n+\frac{1}{2}\) & 0 \cr
        0 & -\De+|x|^2-2\(n+\frac{1}{2}\) \cr
  \end{matrix}
\right).
\end{eqnarray*}
In this case, the representation matrix of $\left.S_{{\mathrm{II}},\om}''(\Om_n)\right|_{\{\Om_{n-2},\Om_{n}, \Om_{n+2}\}}$ is
\begin{eqnarray*}
A_{{\mathrm{II}}}=
\left(
  \begin{matrix}
        -2 & -\frac{\sqrt{n(n-1)}}{4} & 0\cr
        -\frac{\sqrt{n(n-1)}}{4} & -\(n+\frac{1}{2}\) & -\frac{\sqrt{(n+1)(n+2)}}{4} \cr
        0 & -\frac{\sqrt{(n+1)(n+2)}}{4} & 2
  \end{matrix}
\right),
\end{eqnarray*}
and the characteristic polynomial becomes as
\begin{eqnarray*}
F_{\mathrm{II}}(\la)=\la^3+\(n+\frac{1}{2}\)\la^2-\frac{1}{8}\(n^2+n+33\)\la-\frac{7}{2}\(n+\frac{1}{2}\).
\end{eqnarray*}
Since $F_{\mathrm{II}}(\la)=-F_{\mathrm{I}}(-\la)$,
we see that $A_{\mathrm{II}}$ has two negative eigenvalues and one positive eigenvalue.

Therefore, we have the conclusion of the Lemma.
\end{proof}

\subsection*{Acknowledgment}
The authors express their deep gratitude to Professor Kenji Yajima
for fruitful discussions and for his leading the authors to the work of Kato.
The second author is supported by Japan Society for the Promotion of Science (JSPS)
Grant-in-Aid for Research Activity Start-up (22840039).

\def\cprime{$'$}
\providecommand{\bysame}{\leavevmode\hbox to3em{\hrulefill}\thinspace}
\providecommand{\MR}{\relax\ifhmode\unskip\space\fi MR }
\providecommand{\MRhref}[2]{%
  \href{http://www.ams.org/mathscinet-getitem?mr=#1}{#2}
}
\providecommand{\href}[2]{#2}

\end{document}